\documentclass[10pt]{article}
\usepackage{amsmath, amsthm, amssymb, amsfonts, mathrsfs, mathtools}
\usepackage{graphicx}
\usepackage{makeidx}
\usepackage[left=2cm,right=2cm,top=2cm,bottom=2cm]{geometry}
\usepackage{caption}
\usepackage{subcaption}
\usepackage[section]{placeins}
\usepackage{enumitem}
\usepackage{tikz}
\usepackage{stmaryrd}
\usepackage{authblk}
\allowdisplaybreaks

\usepackage{tikz}
\usetikzlibrary{arrows.meta}

\usetikzlibrary{decorations.pathreplacing}
\tikzstyle{vertex}=[auto=left,circle,draw=black,fill=white, inner sep=1.5]

\usepackage{hyperref}
\hypersetup{colorlinks=true, citecolor={blue}, linkcolor=blue, filecolor=magenta, urlcolor=cyan}

\newtheorem{theorem}{Theorem}[section]

\newtheorem{lema}[theorem]{Lemma}
\newtheorem{corollary}{Corollary}[theorem]

\newtheorem{remark}[theorem]{Remark}

\newcommand{\Cay}{\operatorname{Cay}}

\title{A Generalization of Ramanujan's Sum over Finite Rings}

\author[1]{Priya Dhankhar\thanks{Corresponding author}}
\author[2]{Sanjay Kumar Singh}
\affil[ ]{\small{\textsuperscript{1,2}Department of Mathematics, Indian Institute of Science Education and Research Bhopal, India.}}
\affil[ ]{ {\textsuperscript{1}priya22@iiserb.ac.in}
\textsuperscript{2}sanjayks@iiserb.ac.in}

\date{}

\begin{document}
	\maketitle
	
	\vspace{-0.3in}
\begin{center}{\textbf{Abstract}}\end{center}

%Let \(R\) be a finite ring with identity, not necessarily commutative, and let $x \in R$. Define $$[x]_{\ell} := \{ y \in R: (y)_{\ell} = (x)_{\ell} \},$$ 
%where $(x)_{\ell}$ denotes the left ideal of $R$ generated by $x$. We introduce a generalized Ramanujan sum in the setting of finite rings and obtain an explicit formula for it. We show that when $x \in R^{\times}$, the generalized Ramanujan sum gives the eigenvalues of the unitary Cayley graph $\mathrm{Cay}(R,R^{\times})$, and hence provides an effective method for computing the spectrum of these graphs. Surprisingly, Under an 
%independence condition on the maximal left ideals of $(x)_{\ell}$, this formula extends the known formula of classical Ramanujan sums to finite rings. This generalization not only enriches the theory of Ramanujan sums but also provides new tools in spectral graph theory, representation theory, and algebraic number theory.

In this paper, we introduce a generalization of Ramanujan sums over finite rings and derive an explicit formula for these sums. These generalized Ramanujan sums coincide
with the eigenvalues of the unitary Cayley graph
$\mathrm{Cay}(R,R^{\times})$, as follows from Babai's result
on Cayley graph spectra. This provides a useful approach to determining the spectrum of the unitary Cayley graph. Moreover, under a certain independence condition on the maximal left ideals, we show that our formula extends the classical formula for Ramanujan sums to finite rings.

\vspace*{0.3cm}

%\textbf{Mathematics Subject Classifications:} Primary 11L05, 05C50; Secondary 20C15, 05C25.	
	
\section{Introduction}

Ramanujan's sum is a classical arithmetic function over the integers, first introduced by Srinivasa Ramanujan in the early 20th century. It has important applications in number theory and spectral graph theory, particularly in the study of circulant graphs. In this paper, we study a generalization of Ramanujan's sum over finite rings.

Throughout the paper, we consider $(R,+, \times )$ to be a finite ring with unity. We will use the notation $R^{\times}$ to denote the set of units of $R$. For any $x \in R$, we denote by $(x)_{\ell}$ the left ideal of $R$ generated by $x$. Define $$[x]_{\ell} := \{ y \in R: (y)_{\ell} = (x)_{\ell} \}.$$

The set \([x]_{\ell}\) is a natural ring-theoretic analogue of the reduced
residue classes modulo \(n\). Indeed, in \(\mathbb{Z}/n\mathbb{Z}\), elements
generating the same principal ideal are precisely those having the same
greatest common divisor with \(n\). Furthermore, when \(x\in R^{\times}\), we have
\([x]_{\ell}=R^{\times}\). 

Since the additive group of $R$ is isomorphic to $\mathbb{Z}_{n_1} \times \cdots \times \mathbb{Z}_{n_k}$, we will identify the additive group $(R,+)$ with the abelian group $\mathbb{Z}_{n_1} \times \cdots \times \mathbb{Z}_{n_k}$. For any element $\alpha \in R$, we write $\alpha := (\alpha_1,\ldots , \alpha_k)$, where $\alpha_i \in \mathbb{Z}_{n_i}$ for each $i \in \{ 1,\ldots k \}$. Note that each additive character $\psi$ of $R$ can be expressed as $\psi_{\alpha}$ for some $\alpha \in R$ such that 
\begin{equation*}
\psi_{\alpha}(x)=\prod_{j=1}^{k}\omega_{n_j}^{\alpha_j x_j} \textnormal{ for each } x=(x_1,\ldots ,x_k) \in \mathbb{Z}_{n_1} \times \cdots \times \mathbb{Z}_{n_k} \textnormal{ and } \omega_{n_j}=  \exp{\left(\frac{2\pi \textbf{i}}{n_j}\right)}. 
\end{equation*}

%In this paper, we will use the following terminologies:
%\begin{itemize}
%\item Since the additive group of $R$ is isomorphic to $\mathbb{Z}_{n_1} \times \cdots \times \mathbb{Z}_{n_k}$, we will identify the additive group $(R,+)$ with the abelian group $\mathbb{Z}_{n_1} \times \cdots \times \mathbb{Z}_{n_k}$.
%\item For any element $\alpha \in R$, we write $\alpha := (\alpha_1,\ldots , \alpha_k)$, where $\alpha_i \in \mathbb{Z}_{n_i}$ for each $i \in \{ 1,\ldots k \}$.
%\item For any $\alpha \in R$, we denote by $\langle \alpha \rangle$ the cyclic subgroup of the additive group of $R$ generated by $\alpha$.
%\item For any $\alpha \in R$, define the additive character of $R$ $$\psi_{\alpha}: R \to \mathbb{C}^{\times}$$ such that 
%\begin{equation*}
%\psi_{\alpha}(x)=\prod_{j=1}^{k}\omega_{n_j}^{\alpha_j x_j} \textnormal{ for each } x=(x_1,\ldots ,x_k) \in \mathbb{Z}_{n_1} \times \cdots \times \mathbb{Z}_{n_k} \textnormal{ and } \omega_{n_j}=  \exp{(\frac{2\pi \iota}{n_j})}. 
%\end{equation*} Note that all the additive character of $R$ is of the form of $\psi_{\alpha}$ for some $\alpha \in R$.
%\end{itemize}

\noindent For any $\alpha , x \in R$, we define a generalization of Ramanujan's sum over finite Rings as 
\begin{align}
C_{\alpha}(x):= \sum_{s \in [x]_{\ell} } \psi_{\alpha}(s).\nonumber
\end{align}

 Recently, Ramanujan's sum was generalized to finite groups in \cite{RamaGenKadyan}, where elements are grouped according to the normal subgroups generated by their conjugacy classes. In this paper, we introduce a similar generalization for finite rings by grouping together elements that generate the same principal left ideal.
%A related generalization of Ramanujan's sum over finite groups was recently introduced in \cite{RamaGenKadyan}, using equivalence classes determined by normal subgroups generated by conjugacy classes. In contrast, our construction is based on elements generating the same principal left ideal of a finite ring.

\textbf{Classical Ramanujan's sum: } For the special case where $R=\mathbb{Z}_n$ is the ring of integers modulo $n$, equipped with the usual addition and multiplication, and $x\equiv1(\mod n)$, then $$C_{\alpha}(1)= \sum_{ \substack{ 1\leq s \leq n\\ \gcd(s,n)=1} } \omega_n^{s\alpha} = R_n(\alpha).$$
Here, $R_n(\alpha)$ denotes the classical Ramanujan's sum, which was introduced by Ramanujan in 1918. For further details, see \cite{ramanujan1918certain}. It is known that 
\begin{align}\label{neweqRamaNujansum}
   R_n(\alpha)=  \mu\bigg({\frac{n}{\delta_{\alpha}}}\bigg) \frac{\varphi(n)}{\varphi(\frac{n}{\delta_{\alpha}})}, 
\end{align}
where $\delta_{\alpha} = \gcd(n,\alpha)$, $\mu$ is classic M$\ddot{\text{o}}$bius function, and $\varphi$ is Euler's phi function. This naturally leads to the question of whether a similar formula can be derived for $C_{\alpha}(x)$, where $R$ is an arbitrary finite ring. In this paper, we answer this question in the affirmative.  

\textbf{Unitary Cayley graph}: In 1979, Babai \cite{babai1979spectra} used group characters to determine the eigenvalues of Cayley graphs. His result implies that  the eigenvalue of the Cayley graph $\Cay(R , [x]_{\ell})$ is given by
\begin{align}
\lambda_{\alpha}= C_{\alpha}(x)~ \textnormal{ for each } \alpha \in R. \nonumber
\end{align} If $x\in R^{\times}$, then the set $\{ C_{\alpha}(x): \alpha \in R \}$  represents all the eigenvalues of the unitary Cayley graph $\Cay(R , R^{\times})$.

Thus, the character sum \(C_{\alpha}(x)\) has interesting connections with classical Ramanujan sums and the
eigenvalues of unitary Cayley graphs. These connections highlight the importance of studying \(C_{\alpha}(x)\). Our contributions of this paper are as follows.

\begin{enumerate}
    \item \textbf{A formula for arbitrary finite rings:} For \(\alpha,x\in R\), we express \(C_{\alpha}(x)\) using the maximal $R$-left
    ideals of \((x)_{\ell}\) and the largest left ideal contained in \((x)_{\ell}\cap \ker\psi_{\alpha}\). The result is as follows.\\
\noindent \textbf{Theorem} \textit{(Theorem~\ref{MainResultProof})} Let $R$ be a finite ring and $\alpha, x \in R$. Then 
\begin{align}
C_{\alpha}(x) = \left\{ \begin{array}{cl}
		  |[x]_{\ell}|  &  \mbox{if } K=(x)_{\ell}\\ 
	    \sum\limits_{ \substack{ E\subseteq \mathcal{M}_R((x)_{\ell}) \\ \bigcap\limits_{M \in E}M \subseteq  K } }  |\bigcap\limits_{M \in E}M|~  (-1)^{|E|}   & \mbox{if } \bigcap\limits_{M \in \mathcal{M}_R((x)_{\ell}) } M \subseteq K \\
        0 &  \mbox{otherwise}  \\ 
	\end{array}\right. \nonumber
\end{align} 
where $K$ is the largest left ideal of $R$ contained in $(x)_{\ell} \cap \{ y \in R : \psi_\alpha(y)=1 \}$. Here, $\mathcal{M}_R((x)_{\ell})$ to be the set of all largest proper left ideals of  $(x)_{\ell}$ that are also left ideals of $R$. 

    \item \textbf{A Ramanujan-type closed formula:} Under a natural independence condition on the maximal $R$-left ideals of
    \((x)_{\ell}\), we obtain a formula for $C_{\alpha}(x)$. This formula extends the classical
Ramanujan sum to the setting of finite rings. The result is as follows.\\
    \noindent \textbf{Theorem} \textit{(Theorem~\ref{MainResultProof00})} Let $R$ be a finite ring and $\alpha, x \in R$. If there is no proper subset $E$ of $\mathcal{M}_{R}((x)_{\ell})$ with $\bigcap\limits_{M \in E}M = \bigcap\limits_{M \in \mathcal{M}_{R}((x)_{\ell}) }M$, then
\begin{align}
C_{\alpha}(x) =\mu_R( K , (x)_{\ell} ) ~~ \frac{|[x]_{\ell}|}{\varphi_R( K , (x)_{\ell} )}. \nonumber
\end{align} 
where $K$ is the largest left ideal of $R$ contained in $(x)_{\ell} \cap \{ y \in R : \psi_\alpha(y)=1 \}$.  The definitions of the functions $\varphi_R$ and $\mu_R$ are provided in Section 3 and Section 4, respectively.
\end{enumerate}

The remainder of the paper is organized as follows. In Section~2, we review the background and related results that motivate the
study of \(C_{\alpha}(x)\). In Section~3, we define the relevant
Möbius function on the lattice of left ideals and establish the
preliminary results used later. In Section~4, we derive the explicit
formulas for \(C_\alpha(x)\) and obtain the commutative-ring version as
a corollary. In Section~5, we discuss analogous formulations involving right ideals
and two-sided ideals.

%%%%%%%%%%%%%%%%%%%%%%%%%%%%%%%%%%%%%%%%%%%%%%%%%%%%%%%%%%%%%%%%%%%%%%%%%%%%%%%%%%%%%%%%%%%%%%%%%%%%%%%%%%%%%%%%%%%%%%%%%%%%%%
%%%%%%%%%%%%%%%%%%%%%%%%%%%%%%%%%%%%%%%%%%%%%%%%%%%%%%%%%%%%%%%%%%%%%%%%%%%%%%%%%%%%%%%%%%%%%%%%%%%%%%%%%%%%%%%%%%%%%%%%%%%%%%
%%%%%%%%%%%%%%%%%%%%%%%%%%%%%%%%%%%%%%%%%%%%%%%%%%%%%%%%%%%%%%%%%%%%%%%%%%%%%%%%%%%%%%%%%%%%%%%%%%%%%%%%%%%%%%%%%%%%%%%%%%%%%%
%%%%%%%%%%%%%%%%%%%%%%%%%%%%%%%%%%%%%%%%%%%%%%%%%%%%%%%%%%%%%%%%%%%%%%%%%%%%%%%%%%%%%%%%%%%%%%%%%%%%%%%%%%%%%%%%%%%%%%%%%%%%%%

%%%%%%%%%%%%%%%%%%%%%%%%%%%%%%%%%%%%%%%%%%%%%%%%%%%%%%%%%%%%%%%%%%%%%%%%%%%%%%%%%%%%%%%%%%%%%%%%%%%%%%%%%%%%%%%%%%%%%%%%%%%%%%
%%%%%%%%%%%%%%%%%%%%%%%%%%%%%%%%%%%%%%%%%%%%%%%%%%%%%%%%%%%%%%%%%%%%%%%%%%%%%%%%%%%%%%%%%%%%%%%%%%%%%%%%%%%%%%%%%%%%%%%%%%%%%%
%%%%%%%%%%%%%%%%%%%%%%%%%%%%%%%%%%%%%%%%%%%%%%%%%%%%%%%%%%%%%%%%%%%%%%%%%%%%%%%%%%%%%%%%%%%%%%%%%%%%%%%%%%%%%%%%%%%%%%%%%%%%%%
%%%%%%%%%%%%%%%%%%%%%%%%%%%%%%%%%%%%%%%%%%%%%%%%%%%%%%%%%%%%%%%%%%%%%%%%%%%%%%%%%%%%%%%%%%%%%%%%%%%%%%%%%%%%%%%%%%%%%%%%%%%%%%

\section{Preliminaries}

In this section, we introduce some basic definitions and results from ring theory that will be used throughout this article. The following terminology and notations will be used in this paper:
\begin{itemize}
    \item  We will use the notation $\mathcal{I}(R)$ to denote the collection of all left ideals of $R$. 
    \item For any $I \in \mathcal{I}(R)$, define $$[I^R]:= I \setminus \bigcup\limits_{\substack{J \in \mathcal{I}(R)\\ J \subsetneq I }} J.$$ Let $x\in R$. We observe that $[x] = [I^R]$ whenever $I=(x)_{\ell}$. Similarly, if $I$ is not a principal left ideal of $R$ then $[I^R]=\emptyset$.
    
    \item For any $I,J \in \mathcal{I}(R)$ and $I \subseteq J$, we call $I$ to be an \textit{$R$-left ideal} of $J$ if $I$ is a left ideal of $R$.

    \item For any $I,J \in \mathcal{I}(R)$ and $I \subseteq J$, we call $I$ to be a \textit{maximal $R$-left ideal} of $J$ if it is a $R$-left ideal of $J$ and it is the largest $R$-left ideal of $J$ that is properly contained in $J$.
    
   \item For any $J \in \mathcal{I}(R)$, we denote by $\mathcal{M}_R(J)$ the collection of all maximal $R$-left ideals of $J$.
   \item For any $I,J \in \mathcal{I}(R)$, we denote by  $\mathcal{M}_R(J,I)$ the collection of all maximal $R$-left ideals of $J$ which contains $I$.

    \item Let $J \in \mathcal{I}(R)$ and $\mathcal{M}_1, \mathcal{M}_2 \subseteq \mathcal{M}_R(J)$. We say that $\mathcal{M}_1$ is a minimal subset of $\mathcal{M}_2$ if it satisfies the following three conditions:
\begin{enumerate}[label=(\roman*)]
    \item $\mathcal{M}_1 \subseteq  \mathcal{M}_2$,
    \item $\bigcap\limits_{M \in \mathcal{M}_1}M = \bigcap\limits_{M \in \mathcal{M}_2}M,$
    \item $\bigcap\limits_{M \in E}M \neq \bigcap\limits_{M \in \mathcal{M}_1}M \mbox{ for all } E \subsetneq \mathcal{M}_1$.
\end{enumerate}
     If $\mathcal{M}_1$ is a minimal subset of $\mathcal{M}_1$, then we say $\mathcal{M}_1$ is a minimal subset of itself.
\begin{remark}
    If $\mathcal{M}_1$ is a minimal subset of $\mathcal{M}_2$, then $\mathcal{M}_1$ is a minimal subset of itself, and every subset of $\mathcal{M}_1$ is a minimal subset of itself.
\end{remark}

    \item For any $I, J \in \mathcal{I}(R)$ and $I \subseteq J$, define 
$$\varphi_R(I,J):= \left\{ \begin{array}{cl}
		1 & \mbox{if }I=J  \\
	    \prod\limits_{M \in \mathcal{M}_R(J,I)} \left( \frac{|J|}{|M|} - 1 \right)    & \mbox{if }  I =\bigcap\limits_{M\in \mathcal{M}_{R}(J,I)}M_{i} \\
		1   &\mbox{otherwise} 
	\end{array}\right. .$$ 
\end{itemize}

%{ can be express as an intersection  of }\\ & \mbox{ some maximal $R$-left ideals of  } J 

\begin{lema}\label{NewLemaForMainThm}
         Let $R$ be a finite ring, $I$ be a left ideal of $R$, and $K$ can be express as an intersection of some maximal $R$-left ideals of $I$. Then $$\varphi_R(K,I) = \sum_{E \subseteq \mathcal{M}_R(I,K)} (-1)^{|\mathcal{M}_R(I,K)|-|E|} \prod_{M \in E} \frac{|I|}{|M|}.$$ 
\end{lema}
\begin{proof} We have
    \begin{align*}
     \varphi_R(K,I) &= \prod_{M \in \mathcal{M}_R(I,K)} \left( \frac{|I|}{|M|} - 1 \right) \\ 
     &= \sum_{E \subseteq \mathcal{M}_R(I,K)} (-1)^{|\mathcal{M}_R(I,K)|-|E|} \prod_{M \in E} \frac{|I|}{|M|}. 
    \end{align*}
   
\end{proof}

%%%%%%%%%%%%%%%%%%%%%%%%%%%%%%%%%%%%%%%%%%%%%%%%%%%%%%%%%%%
%%%%%%%%%%%%%%%%%%%%%%%%%%%%%%%%%%%%%%%%%%%%%%%%%%%%%%%%%%%%
In general, the Chinese Remainder Theorem holds for collections of two-sided ideals. However, it does not necessarily hold for left ideals. For example, consider the ring $R = M_2(\mathbb{Z}_2)$, the ring of $2\times 2$ matrices over the field $\mathbb{Z}_2$. Define three subsets of $R$ as follows:
\begin{itemize}
    \item Let $I_1$ be the set of matrices whose rows lie in $\{(1,0),(0,0)\}$.
    \item Let $I_2$ be the set of matrices whose rows lie in $\{(0,1),(0,0)\}$.
    \item Let $I_3$ be the set of matrices whose rows lie in $\{(1,1),(0,0)\}$.
\end{itemize} These are maximal left ideals of $R$. However, one can verify that $$\frac{|R|}{ |I_1 \cap I_2 \cap I_3|} \neq \frac{|R|}{ |I_1|} \frac{|R|}{ |I_2|} \frac{|R|}{ |I_3|}.$$  Thus, the classical version of the Chinese Remainder Theorem fails in this case. However, a variant of the theorem can hold for left ideals if an additional condition of minimality is imposed on the collection of maximal $R$-left ideals. The following result provides an analogue of the Chinese Remainder Theorem in the context of left ideals.

\begin{lema}\label{LemmaMaximEquality00}
    Let $R$ be a finite ring, $I$ be a left ideal of $R$, and $\mathcal{M}_1, \mathcal{M}_2$ be disjoint subsets of  $\mathcal{M}_R(I)$.
    \begin{enumerate}[label=(\roman*)]
        \item $\mathcal{M}_1$ is minimal subset of itself if and only if 
        \begin{align}
        \frac{|I|}{|\bigcap\limits_{M \in \mathcal{M}_1} M|} = \prod\limits_{M \in \mathcal{M}_1} \frac{|I|}{|M|}. \label{CRTeqLemma}
        \end{align}
        \item If $\mathcal{M}_R(I)$ is minimal subset of itself then $$\frac{|I|}{|\bigcap\limits_{M \in E_1 \cup E_2 } M|} = \frac{|I|}{|\bigcap\limits_{M \in E_1} M|} ~ \frac{|I|}{|\bigcap\limits_{M \in E_2} M|}$$  for all $ E_1 \subseteq \mathcal{M}_1$  and $ E_2 \subseteq \mathcal{M}_2$.
    \end{enumerate}
\end{lema}
\begin{proof}
\begin{enumerate}[label=(\roman*)]
    \item Assume that $\mathcal{M}_1$ is a minimal subset of itself, and let $|\mathcal{M}_1 |= k$. We proceed by induction on $k$. If $k=1$, then the result holds trivially. Assume, as the induction hypothesis, that the result is true for $k=t-1$. Now consider the case $k=t$. Let $K$ be a fix element of $\mathcal{M}_1$, and define $\mathcal{M}_2 := \mathcal{M}_1 \setminus \{K\}$.  Since $\mathcal{M}_1$ is minimal subset of itself, it follows that the sum of the ideals $\bigcap\limits_{M \in \mathcal{M}_2}M$ and $K$ is equal to $I$. The Second Isomorphism Theorem for rings, we have
        \begin{align}\label{IntersecEqua1}
        \frac{|I|}{|\bigcap\limits_{M \in \mathcal{M}_2} M|} = \frac{|K|}{|\bigcap\limits_{M \in \mathcal{M}_1} M|}.
        \end{align}
        Now, applying the induction hypothesis to $\mathcal{M}_2$, Equation~(\ref{IntersecEqua1}) yields 
        \begin{align}\label{IntersecEqua2}
        \prod\limits_{M \in \mathcal{M}_2} \frac{|I|}{|M|}  = \frac{|K|}{|\bigcap\limits_{M \in \mathcal{M}_1} M|}.
        \end{align}
        Multiplying both sides of Equation~(\ref{IntersecEqua2}) by $\frac{|I|}{|K|}$, we obtain the desired result. Conversely, assume tha Equation (\ref{CRTeqLemma}) holds, and suppose that $\mathcal{M}_1$ is not minimal subset of itself. Then, there exists a proper subset $\mathcal{M}_2$ of $\mathcal{M}_1$ that is minimal in $\mathcal{M}_1$, hence also minimal in itself. As $\mathcal{M}_2$ is minimal in itself, therefore from the above proof we have 
        \begin{align} \label{CRTforSubsetM2}
            \frac{|I|}{|\bigcap\limits_{M\in \mathcal{M}_2}|}=\prod \limits_{M\in \mathcal{M}_2}\frac{|I|}{|M|}.
        \end{align}     
Note that $\mathcal{M}_2$ is minimal subset of $\mathcal{M}_1$, therefore $\bigcap\limits_{M\in \mathcal{M}_2} M= \bigcap\limits_{M\in \mathcal{M}_1} M$. Thus,
$$ \frac{|I|}{\big | \bigcap\limits_{M\in \mathcal{M}_2} M \big |} = \frac{|I|}{\big | \bigcap\limits_{M\in \mathcal{M}_1} M \big |}.$$
Upon comparing Equations (\ref{CRTeqLemma}) and (\ref{CRTforSubsetM2}), we get
$$ \prod \limits_{M\in \mathcal{M}_1} \frac{|I|}{|M|}= \prod \limits_{M\in \mathcal{M}_2} \frac{|I|}{|M|}. $$ Which implies 
$ \prod \limits_{M\in \mathcal{M}_1 \setminus \mathcal{M}_2} \frac{|I|}{|M|} =1$, therefore $|I| =|M|$ for all $M \in \mathcal{M}_1\setminus M_2$. However, this is not possible because each $M$ is a proper subset of $I$. Therefore, our assumption is wrong, and we conclude that $\mathcal{M}_1$ is minimal in itself.
        
    \item  $\mathcal{M}_R(I)$ is minimal in itself, implies that $E_1$ and $E_2$ are minimal in $\mathcal{M}_R(I)$. Since $E_1$ and $E_2$ are disjoint subset, their $E_1 \cup E_2$ is also minimal subset in $\mathcal{M}_R(I)$, hence minimal in itself. Therefore, the result follows from part $(i)$.
    
    \end{enumerate}
\end{proof}

%%%%%%%%%%%%%%%%%%%%%%%%%%%%%%%%%%%%%%%%%%%%%%%%%%%%%%%%%%%%%%%%%%%%%%%%%%%%%%%%%%%%%%%%%%%%%%%%%%%%%%%%%%%%%%%%%%%%%%%%%%%%%%%%%%%%%%%%%%%%%%%%%%%%%%%%%%%%%%%%%%%%%%%%%%%%%%%%%%%%%%%%%%%%%%%%%%%%%%%%%%%%%%%%%%%%%%%%%%%%%%%%%%%%%%%%%%%%%%%%%%%%%%%%%%%%%%%%%%%%%%%%%%%%%%%%%%%%%%%%%%%%%%%%%%%%%%%%%%%%%%%%%%%%%%%%%%%%%%%%%%%%%%%%%%%%%%%%
\begin{lema}\label{LemaOnMaxIdealOfCapMaxiIdeal}
    Let $R$ be a ring, $I$ be a left ideal of $R$, $\mathcal{M}_1 \subseteq \mathcal{M}_R(I)$ be a subset that is minimal in itself, and $K \in \mathcal{M}_R(I)$. If $\bigcap \limits_{M\in \mathcal{M}_1 \cup \{K\}}M $ is properly contained in $\bigcap \limits_{M\in \mathcal{M}_1}M$, then  $\bigcap \limits_{M\in \mathcal{M}_1 \cup \{K\}}M$ is maximal $R$-left ideal of $\bigcap \limits_{M\in \mathcal{M}_1}M$.
\end{lema}

\begin{proof}
Let \( L \) be a \( R \)-left ideal such that
\[
\bigcap_{M \in \mathcal{M}_1 \cup \{K\}} M \subsetneq L \subseteq \bigcap_{M \in \mathcal{M}_1} M.
\]
We may assume that \( L \nsubseteq K \) and \( \bigcap\limits_{M \in \mathcal{M}_1} M \nsubseteq \{K\} \). Otherwise:
\begin{enumerate}[label=(\roman*)]
    \item If \( L \subseteq K \), then \( L \subseteq \bigcap \limits_{M \in \mathcal{M}_1\cup \{K\}} M \), contradicting our assumption that \( L \) properly contains this intersection.
    \item If \( \bigcap \limits_{M \in \mathcal{M}_1} M \subseteq K \), then \( \bigcap\limits_{M \in \mathcal{M}_1 \cup \{K\}} M = \bigcap \limits_{M \in \mathcal{M}_1} M \), again contradicting the assumption that the intersection is proper.
\end{enumerate}
 
\noindent From the Second Isomorphism Theorem of rings, we have
$$ \frac{|L+K|}{|K|}=\frac{|L|}{|L \cap K|}.$$
Since $K$ is a maximal $R$-left ideal of $I$ and $L \nsubseteq K$ it follows that $L+K=I$. Therefore, we get
    \begin{align} \label{ord_of_I_from_L}
        |I|=\frac{|L||K|}{|L\cap K|}
    \end{align}

\noindent Again apply the Second Isomorphism Theorem of rings, we have
$$ \frac{\big |\bigcap \limits_{M\in \mathcal{M}_1}M +K \big |}{|k|}=\frac{\big |\bigcap \limits_{M\in \mathcal{M}_1}M \big |}{ \big | \bigcap \limits_{M\in \mathcal{M}_1 \cup \{K\}}M \big |}.$$
Since $\bigcap \limits_{M\in \mathcal{M}_1}M \not\subseteq  K$ and $K$ is a maximal $R$-left ideal of $I$, it follows that $ \bigcap \limits_{M\in \mathcal{M}_1}M +K = I$. Which implies
\begin{align} \label{ord_of_I_from_maxiaml_ideal}
    |I|=\frac{\big |\bigcap \limits_{M\in \mathcal{M}_1}M \big |\big |K\big |}{ \big | \bigcap \limits_{M\in \mathcal{M}_1}M \cap K\big |}.
\end{align}
Upon comparing Equation (\ref{ord_of_I_from_L}) and  Equation (\ref{ord_of_I_from_maxiaml_ideal}), we get 
$$\frac{|L||K|}{|L\cap K|}= \frac{\big |\bigcap \limits_{M\in \mathcal{M}_1}M \big |\big |K\big |}{ \big | \bigcap \limits_{M\in \mathcal{M}_1\cup \{K\}}M \big |}.$$
Therefore
\begin{align} \label{equ_of_order_of_L_and_lcapK}
    \frac{|L|}{\big |\bigcap \limits_{M\in \mathcal{M}_1}M\big |}= \frac{|L \cap K|}{\big |\bigcap \limits_{M\in \mathcal{M}_1\cup \{K\}}M \big |}.
\end{align}
Since $L \subseteq \bigcap \limits_{M\in \mathcal{M}_1}M $, it follows that $L \cap K \subseteq \bigcap \limits_{M\in \mathcal{M}_1 \cup \{K\}}M$, and also $\bigcap \limits_{M\in \mathcal{M}_1\cup \{K\}}M \subseteq L$ which implies $$\bigcap \limits_{M\in \mathcal{M}_1 \cup \{K\}}M \subseteq L\cap K,$$ and hence $$L \cap K = \bigcap \limits_{M\in \mathcal{M}_1 \cup \{K\}}M.$$
It follows that right hand side of Equation (\ref{equ_of_order_of_L_and_lcapK}) is $1$, this forces $|L|=|\bigcap \limits_{M\in \mathcal{M}_1}M|$. This implies $L=\bigcap \limits_{M\in \mathcal{M}_1}M$. Therefore, $\bigcap \limits_{M\in \mathcal{M}_1 \cup \{K\}}M $ is maximal $R$-left ideal of $\bigcap \limits_{M\in \mathcal{M}_1}M$.

\end{proof}

%%%%%%%%%%%%%%%%%%%%%%%%%%%%%%%%%%%%%%%%%%%%%%%%%%%%%%%%%%%%%%%%%%%%%%%%%%%%%%%%%%%%%%%%%%%%%%%%%%%%%%%%%%%%%%%%%%%%%%%%%%%%%%%%%%%%%%%%%%%%%%%%%%%%%%%%%%%%%%%%%%%%%%%%%%%%%%%%%%%%%%%%%%%%%%%%%%%%%%%%%%%%%%%%%%%%%%%%%%%%%%%%%%%%%%%%%%%%%%%%%%%%%%%%%%%%%%%%%%%%%%%%%%%%%%%%%%%%%%%%%%%%%%%%%%%%%%%%%%%%%%%%%%%%%%%%%%%%%%%%%%%%%%%%%%%%%%%%
\begin{lema}\label{NewLemEqua22}
    Let $R$ be a finite ring, and let $I,K$ are left ideals of $R$ with $K \subsetneq I$. If $\bigcap\limits_{M \in \mathcal{M}_R(I)} M \subseteq K$, then $K$ can be expressed as an intersection of some maximal $R$-left ideals of $I$. Moreover, $$K=\bigcap_{M \in \mathcal{M}_R(I,K)} M$$. 
\end{lema}
\begin{proof}
Assume that $\mathcal{M}_1$ is minimal subset of $\mathcal{M}_R(I,K)$.  If $$\bigcap\limits_{M \in \mathcal{M}_1} M = \bigcap\limits_{M \in \mathcal{M}_R(I)} M$$ then $K=\bigcap\limits_{M \in \mathcal{M}_1} M$. This is  because $K \subseteq \bigcap\limits_{M \in \mathcal{M}_1} M =\bigcap\limits_{M \in \mathcal{M}_R(I)} M \subseteq K$. Assume that $$\bigcap\limits_{M \in \mathcal{M}_1} M \neq \bigcap\limits_{M \in \mathcal{M}_R(I)} M.$$ We can choose $\mathcal{M}_2$ a subset of $\mathcal{M}_R(I) \setminus \mathcal{M}_R(I,K)$ such that $\mathcal{M}_1 \cup \mathcal{M}_2$ is a minimal subset of $\mathcal{M}_R(I)$.
Note that 
\begin{align}
   \bigcap \limits_{M\in \mathcal{M}_1 \cup \mathcal{M}_2}M \subsetneq K \subseteq \bigcap \limits_{M \in \mathcal{M}_1}M. \label{NeqEqLstLemma1}
\end{align}
 We prove the result by applying induction to the size of $\mathcal{M}_2$. If the size of $\mathcal{M}_2$ is $1$, then by Lemma \ref{LemaOnMaxIdealOfCapMaxiIdeal} $\bigcap\limits_{M \in \mathcal{M}_1 \cup \mathcal{M}_2} M$ is maximal $R$-left ideal of $\bigcap\limits_{M \in \mathcal{M}_1} M$, hence it follows that $K=\bigcap\limits_{M \in \mathcal{M}_1} M$. Assume that the size of $\mathcal{M}_2$ is $t$ and $\mathcal{M}_2= \{ M_1,M_2, \ldots , M_t\}$. Define $K_1:=K \cap M_t$, $\mathcal{M}_1':=\mathcal{M}_1 \cup \{M_t\}$, and $\mathcal{M}_2':=\{ M_1,M_2,\ldots , M_{t-1} \}$. Equation (\ref{NeqEqLstLemma1}) implies that
 \begin{align}
   \bigcap \limits_{M\in \mathcal{M}_1' \cup \mathcal{M}_2'}M \subsetneq K_1 \subseteq \bigcap \limits_{M \in \mathcal{M}_1'}M. \label{NeqEqLstLemma2}
\end{align} Apply induction hypothesis on $K_1$, $\mathcal{M}_1'$, and  $\mathcal{M}_2'$. We get 
$$K_1 = \bigcap\limits_{M \in \mathcal{M}_1'} M.$$ This implies that $$\bigcap\limits_{M \in \mathcal{M}_1'} M =K_1 \subsetneq K \subseteq \bigcap\limits_{M \in \mathcal{M}_1} M.$$ By Lemma \ref{LemaOnMaxIdealOfCapMaxiIdeal}, $\bigcap\limits_{M \in \mathcal{M}_1' } M$ is maximal $R$-left ideal of $\bigcap\limits_{M \in \mathcal{M}_1} M$. Hence $K = \bigcap\limits_{M \in \mathcal{M}_1} M$.
\end{proof}

%%%%%%%%%%%%%%%%%%%%%%%%%%%%%%%%%%%%%%%%%%%%%%%%%%%%%%%%%%%%%%%%%%%%%%%%%%%%%%%%%%%%%%%%%%%%%%%%%%%%%%%%%%%%%%%%%%%%%%%%%%%%%%%%%%%%%%%%%%%%%%%%%%%%%%%%%%%%%%%%%%%%%%%%%%%%%%%%%%%%%%%%%%%%%%%%%%%%%%%%%%%%%%%%%%%%%%%%%%%%%%%%
%%%%%%%%%%%%%%%%%%%%%%%%%%%%%%%%%%%%%%%%%%%%%%%%%%%%%%%%%%%%%%%%%%%%%%%%%%%%%%%%%%%%%%%%%%%%%%%%%%%%%%%%%%%%%%%%
%%%%%%%%%%%%%%%%%%%%%%%%%%%%%%%%%%%%%%%%%%%%%%%%%%%%%%%%%%%%%%%%%%%%%%%%%%%%%%%%%%%%%
%%%%%%%%%%%%%%%%%%%%%%%%%%%%%%%%%%%%%%%%%%%%%%%%%%%%%%%%%%%%%%%%%%%%%%%%%%%%%%%%%%%%%

\section{M$\ddot{\text{o}}$bius Inversion Formula}

%In 1831, August Ferdinand M$\ddot{\text{o}}$bius introduced the classical M$\ddot{\text{o}}$bius inversion formula, which establishes a relationship between two arithmetic functions, where each is defined in terms of the other by taking sums over divisors. Later, in 1964, Gian-Carlo Rota~\cite{rota1964foundations} proposed a generalization of this formula, extending it to functions defined over partially ordered sets. In earlier works, Weisner \cite{weisner1935abstract} and Hall \cite{hall1936eulerian} explored a version of the M$\ddot{\text{o}}$bius function defined on the lattice of subgroups of a group. In this paper, we consider a version of the M$\ddot{\text{o}}$bius function defined on the left ideals of a ring $R$, inspired by Hall’s approach \cite{hall1936eulerian}. We define the function $\mu_R : \mathcal{I}(R) \times \mathcal{I}(R) \to \mathbb{Z}$ as follows:

In 1964, Gian-Carlo Rota~\cite{rota1964foundations} introduced the concept of the incidence algebra of a poset to study the M$\ddot{\text{o}}$bius function and the M$\ddot{\text{o}}$bius inversion formula. He showed that the M$\ddot{\text{o}}$bius function is the unique function that satisfies a certain recursive relation. Rota’s work brought together many ideas in combinatorics and became a foundation for modern algebraic and enumerative combinatorics. Detailed discussions of his work can be found in the book \cite{RichardStanley}. 

In this section, we focus on the M$\ddot{\text{o}}$bius function and the M$\ddot{\text{o}}$bius inversion formula in the context of the inclusion order in the set $\mathcal{I}(R)$. While Rota’s approach provides only a recursive formula for the M$\ddot{\text{o}}$bius function, we present a more explicit definition based on Hall’s approach \cite{hall1936eulerian}. In earlier works, Weisner \cite{weisner1935abstract} and Hall \cite{hall1936eulerian} explored a version of the M$\ddot{\text{o}}$bius function defined on the lattice of subgroups of a group. 

We define the function $\mu_R : \mathcal{I}(R) \times \mathcal{I}(R) \to \mathbb{Z}$ as follows:
$$\mu_R(I,J):= \left\{ \begin{array}{cl}
		1 & \mbox{if }I=J  \\
	    e-o    & \mbox{if }  I  \mbox{ can be express as an intersection  of some}\\ & \mbox{ maximal $R$-left ideal of  } J  \\
		0   &\mbox{otherwise} 
	\end{array}\right. $$ Here, $e$ (resp. $o$) denotes the number of ways to express $I$ as the intersection of an even (respectively, odd) number of maximal $R$-left-ideals of $J$.
\begin{remark} \label{mobius function for minimal set}
    \begin{enumerate}[label=(\roman*)]
        \item  If $I$ can be express as an intersection  of some maximal $R$-left ideal of $J$, then 
        \begin{align*}
            \mu_R(I,J)= \sum_{\substack{ E \subseteq \mathcal{M}_{R}(J,I) \\ I= \bigcap\limits_{M \in E}M }}(-1)^{|E|}.
        \end{align*}
        \item  If $\mathcal{M}_R(J)$ is minimal subset of itself then 
$$\mu_R(I,J):= \left\{ \begin{array}{cl}
		1 & \mbox{if }I=J  \\
	    (-1)^r    & \mbox{if }  I  \mbox{ can be express as an intersection  of some}\\ & \mbox{ maximal $R$-left ideal of  } J  \\
		0   &\mbox{otherwise} 
	\end{array}\right. $$  where $M_1 , M_2 , \ldots, M_r$ are maximal $R$-left ideals of $J$ and $I=M_1 \cap \ldots \cap M_r$. 
    \end{enumerate}
\end{remark}

\begin{lema}\label{muSumSubset10} Let $R$ be a finite ring. Then
$$ \sum_{\substack{ I \in \mathcal{I}( R ) \\ K \subseteq  I  \subseteq J }}  \mu_R(I,J)= \left\{ \begin{array}{cl}
		1 &  \mbox{if } K = J \\
		0   &\mbox{otherwise.} 
	\end{array}\right. $$ 	
	
\end{lema}

\begin{proof}
If $K = J$ then $$ \sum_{\substack{ I\in \mathcal{I}( R ) \\  K \subseteq  I  \subseteq J }}  \mu_R(I,J)=   \mu_R(J,J) = 1.$$ Assume that $K \subsetneq J$.  We have
$$\sum_{\substack{ I\in \mathcal{I}( R ) \\  K \subseteq  I  \subseteq J }}  \mu_R(I,J)= \sum_{\substack{ I\in \mathcal{I}( R ) \\  K \subseteq  I  \subseteq J \\ \mu_R(I,J)\neq 0 }}  \mu_R(I,J).$$ Let $M_1, M_2, \ldots , M_r$ are the only maximal $R$-left ideals of $J$ containing $K$. For any $0 \leq j \leq r$, there are ${r \choose j}$ many intersection of maximal $R$-left ideals exist in $\mathcal{I}( R )$. Therefore, we obtain
$$ \sum_{\substack{ I\in \mathcal{I}( R ) \\ K \subseteq  I  \subseteq J \\ \mu_R(I,J)\neq 0 }}  \mu_R(I,J)  = \sum_{j=0}^{r} {r \choose j} (-1)^j=(1-1)^{r}=0.$$
\end{proof}

The next result is a M$\ddot{\text{o}}$bius inversion formula. It can be seen as a special case of the well known M$\ddot{\text{o}}$bius inversion formula of Gian-Carlo Rota~\cite{rota1964foundations}, by considering the inclusion order on the set $\mathcal{I}(R)$.

\begin{theorem}\label{MobiusInvForm}
Let $R$ be a finite ring and $f,g: \mathcal{I}(R) \to \mathbb{Z}$. If 
\begin{align}
f(J)= \sum_{\substack{ I\in \mathcal{I}(R) \\ I \subseteq J }} g(I) \label{MobiusInvForEq1}
\end{align}
then
\begin{align}
  g(J) = \sum_{\substack{ I\in \mathcal{I}(R) \\ I \subseteq J }} f(I)~ \mu_R(I,J). \label{MobiusInvForEq2}
 \end{align}
\end{theorem}
\begin{proof}
Assume that Equation (\ref{MobiusInvForEq1}) holds. We have
\begin{align}
 \sum_{\substack{ I\in \mathcal{I}(R) \\ I \subseteq J }} f(I) ~ \mu_R(I,J) &=   \sum_{\substack{ I \in \mathcal{I}( R )\\ I \subseteq J }}   \sum_{\substack{ K \in \mathcal{I}( R ) \\ K \subseteq I }} g(K) ~ \mu_R(I,J) \nonumber\\
 &=  \sum_{\substack{ K \in \mathcal{I}( R )\\ K \subseteq J }}  \sum_{ \substack{ I \in \mathcal{I}( R ) \\ K  \subseteq I \subseteq J } }    g(K)  ~ \mu_R(I,J) \nonumber\\
 &=  \sum_{\substack{ K \in \mathcal{I}( R )\\ K \subseteq J }}  g(K) \sum_{  \substack{ I \in \mathcal{I}( R ) \\ K  \subseteq I \subseteq J }  }      \mu_R(I,J) \nonumber\\
 &= g(J).\nonumber
\end{align}
Here the first equality follows from Equation (\ref{MobiusInvForEq1}) and last equality follows from Part $(i)$ of Lemma~\ref{muSumSubset10}.    
\end{proof}

\begin{lema}\label{GroupSizeEqEquiLema} Let $R$ be a finite ring and $ x \in R$.
\begin{align*}
|[x]| = \sum_{ \substack{ I\in \mathcal{I}(R) \\ I \subseteq (x)_{\ell} } }     |I | ~ \mu_R(I ,(x)_{\ell}).   
\end{align*} 
\end{lema}

\begin{proof}
Let $ x \in R$. Since the set $ \bigcup \limits_{\substack{ I \in \mathcal{I}(R) \\ I \subseteq (x)_{\ell} }} [I^{R}]$ is a disjoint union and 
$(x)_\ell= \bigcup \limits_{\substack{ I \in \mathcal{I}(R) \\ I \subseteq (x)_{\ell} }} [I^{R}]$, it follows that
\begin{align*}
|(x)_\ell| = \sum_{\substack{ I \in \mathcal{I}(R) \\ I \subseteq (x)_{\ell} }} |[I^{R}]|. \nonumber
\end{align*} 
Moreover, $|[x]| = |[(x)_{\ell}^{R}]|$. Now, the proof follows from Theorem~\ref{MobiusInvForm}.
\end{proof}

\begin{lema}\label{MuZeroCond}
Let $R$ be a finite ring, $I,J$ are left ideals of $R$ and $I \subsetneq J$. If $I$ cannot be expressed as an intersection of some maximal R-left ideals of $J$  then $\mu_R(K,J)=0$ for all $K \in \mathcal{I}(R)$ and $K \subseteq I$.
\end{lema}
\begin{proof}
Assume that $I$ cannot be expressed as an intersection of any maximal R-left ideals of $J$. Then we have $\mu_R(I,J)=0$. Suppose $\mu_R(K,J)\neq 0$ for some $K \in \mathcal{I}(R)$ and $K \subseteq I$. Then either $K=J$ or $K$ can be expressed as an intersection of some maximal $R$-left ideal of $J$. If $K=J$, then it would imply $I=J$, which is not possible by assumption. On the other hand, by Lemma~\ref{NewLemEqua22}, it follows that $I$ can be expressed as an intersection of some maximal $R$-left ideal of $J$. This leads to a contradiction.
\end{proof}

%%%%%%%%%%%%%%%%%%%%%%%%%%%%%%%%%%%%%%%%%%%%%%%%%%%%%%%%%%%%%%%
%%%%%%%%%%%%%%%%%%%%%%%%%%%%%%%%%%%%%%%%%%%%%%%%%%%%%%%%%%%%%%%
%%%%%%%%%%%%%%%%%%%%%%%%%%%%%%%%%%%%%%%%%%%%%%%%%%%%%%%%%%%%%%%

\section{Main Theorem}
In this section, we present the proofs of the main theorems of this paper. Let $S$ be a subgroup of the additive group of the ring $R$. We define $$S^{\perp} = \{ x: x \in R \mbox{ and } \psi_s(x)=1 \mbox{ for all } s\in S\}.$$ By Theorem 17.3 of~\cite{MartinGordon}, if $S$ is an additive subgroup of $R$, then $|S^{\perp}| = \frac{|R|}{|S|}$.

For example, if $Z=\mathbb{Z}_5 \otimes \mathbb{Z}_5 \otimes \mathbb{Z}_{25}$ then $$\langle (1,0,5) \rangle ^{\perp} = \langle (1 , 1, 4) \rangle \cup \langle (1 , 2, 4) \rangle \cup \langle (1 , 3, 4) \rangle \cup \langle (1 , 4, 4) \rangle  \cup \langle (1 , 0, 4) \rangle \cup \langle (0 ,1, 0) \rangle.$$

\noindent For any $\alpha \in R$ and $I \in \mathcal{I}(R)$, define $$f_{\alpha}(I) := \sum\limits_{s\in I } \psi_{\alpha}(s) \mbox{ and } g_{\alpha}(I) :=  \sum\limits_{s\in [ I^{R} ] } \psi_{\alpha}(s).$$

\noindent Note that
\begin{align} \label{relation of c_alpha_xand g_aplha_x}
    g_{\alpha}( (x)_{\ell} ) = C_{\alpha}(x)\hspace{2mm} \text{for all}  \hspace{2mm}x \in R
\end{align}

\noindent We have
\begin{align}
f_{\alpha}(I) = \left\{ \begin{array}{cl}
		|I| &  \mbox{if }  \alpha \in  I^{\perp} \\
		0   &\mbox{otherwise} 
	\end{array}\right.  \nonumber
\end{align} for all $I \in \mathcal{I}(R)$.
We can also write 
\begin{align}
f_{\alpha}((x)_{\ell}) = \sum\limits_{s\in (x)_{\ell} } \psi_{\alpha}(s)= \sum_{\substack{ I \in \mathcal{I}(R) \\ I \subseteq (x)_{\ell} } }   \sum\limits_{s\in [ I^{R} ] } \psi_{\alpha}(s)=  \sum_{\substack{ I \in \mathcal{I}(R) \\ I \subseteq (x)_{\ell} } } g_{\alpha}(I). \nonumber
\end{align}
Here, the second equality follows from the fact that $(x)_\ell= \bigcup \limits_{\substack{ I \in \mathcal{I}(R) \\ I \subseteq (x)_{\ell} }} [I^{R}]$.

\noindent By Theorem \ref{MobiusInvForm}, we get
\begin{align}
g_{\alpha}((x)_{\ell}) = \sum_{ \substack{ I \in \mathcal{I}(R) \\ I \subseteq (x)_{\ell} } } f_{\alpha}(I) ~ \mu_R(I,(x)_{\ell} )
= \sum_{ \substack{ I \in \mathcal{I}(R) \\ I \subseteq (x)_{\ell} \\ \alpha \in I^{\perp} } }  |I|~  \mu_R( I , (x)_{\ell} ) 
&= \sum_{ \substack{ I \in \mathcal{I}(R) \\ I \subseteq (x)_{\ell} \\ I \subseteq \langle \alpha \rangle^{\perp} } }  |I|~  \mu_R( I , (x)_{\ell} )\nonumber\\
&= \sum_{ \substack{ I \in \mathcal{I}(R) \\ I \subseteq (x)_{\ell} \cap \langle \alpha \rangle^{\perp}  } }  |I|~  \mu_R( I , (x)_{\ell} ).\nonumber
\end{align}

\noindent The third equality follows from the fact that $\alpha \in  I^{\perp}$ if and only if $I \subseteq \langle \alpha \rangle^{\perp}$. From Equation (\ref{relation of c_alpha_xand g_aplha_x}), we obtain
\begin{align}
C_{\alpha}( x ) &= \sum_{ \substack{ I \in \mathcal{I}(R) \\ I \subseteq (x)_{\ell} \cap \langle \alpha \rangle^{\perp}  } }  |I|~  \mu_R( I , (x)_{\ell} ).\label{SubFormulaEq1}
\end{align}

Note that the set $(x)_{\ell} \cap \langle \alpha \rangle^{\perp}$ is a subgroup of the additive group of $R$, but it may or may not be a left ideal of $R$.

Define
\[
K
:=
\sum_{\substack{I\in\mathcal{I}(R)\\
I\subseteq (x)_{\ell}\cap \langle \alpha\rangle^{\perp}}} I.
\]
Since $(x)_{\ell}\cap \langle \alpha\rangle^{\perp}$ is an additive subgroup of $R$, the sum of left ideals contained
in $(x)_{\ell}\cap \langle \alpha\rangle^{\perp}$ is again contained in $(x)_{\ell}\cap \langle \alpha\rangle^{\perp}$. Hence $K$ is a left ideal of $R$ contained
in $(x)_{\ell}\cap \langle \alpha\rangle^{\perp}$. Moreover, every left ideal $I$ of $R$ satisfying $I\subseteq (x)_{\ell}\cap \langle \alpha\rangle^{\perp}$ is
contained in $K$. Therefore, $K$ is the largest left ideal of $R$
contained in $(x)_{\ell}\cap \langle \alpha\rangle^{\perp}$.

%%%%%%%%%%%%%%%%%%%%%%%%%%%%%%%%%%%%%%%%%%%%%%%%%%%%%%%%%%%%
%%%%%%%%%%%%%%%%%%%%%%%%%%%%%%%%%%%%%%%%%%%%%%%%%%%%%%%%%%%%
\noindent Equation~(\ref{SubFormulaEq1}) implies that
\begin{align}
C_{\alpha}( x ) = \sum_{ \substack{ I \in \mathcal{I}(R) \\ I \subseteq  K } }  |I|~  \mu_R( I , (x)_{\ell} ), \label{NewEqUniqeKaddeq}
\end{align}
where $K$ is the largest left ideal of $R$ contained in $(x)_{\ell} \cap \langle \alpha \rangle^{\perp}$. 

\begin{lema}\label{NewLemmaqqq} Let $R$ be a finite ring and $ \alpha, x  \in R$.  Then
\begin{align}
C_{\alpha}(x) =  \sum_{ \substack{ E\subseteq \mathcal{M}_R((x)_{\ell}) \\ \bigcap\limits_{M \in E}M \subseteq  K } }  |\bigcap\limits_{M \in E}M|~~  (-1)^{|E|}  \nonumber
\end{align} 
where $K$ is the largest left ideal of $R$ contained in $(x)_{\ell} \cap \langle \alpha \rangle^{\perp}$. 
\end{lema}
\begin{proof} Assume that $K$ is the largest left ideal of $R$ contained in $(x)_{\ell} \cap \langle \alpha \rangle^{\perp}$. 
From Equation~(\ref{NewEqUniqeKaddeq}), we have
\begin{align}
C_{\alpha}( x ) &= \sum_{ \substack{ I \in \mathcal{I}(R) \\ I \subseteq  K \\ \mu_R( I , (x)_{\ell} ) \neq 0 } }  |I|~  \mu_R( I , (x)_{\ell} )\nonumber\\
&= \sum_{ \substack{ I \in \mathcal{I}(R) \\ I \subseteq  K } }  |I|~  \sum\limits_{\substack{ E\subseteq \mathcal{M}_R((x)_{\ell}) \\  I =\bigcap\limits_{M \in E}M }} (-1)^{|E|} \nonumber\\
& = \sum_{ \substack{ E\subseteq \mathcal{M}_R((x)_{\ell}) \\ \bigcap\limits_{M \in E}M \subseteq  K } }  |\bigcap\limits_{M \in E}M|~~  (-1)^{|E|}.\nonumber
\end{align}
Here the second equality follows from Part $(i)$  of Remark \ref{mobius function for minimal set}.
\end{proof}

\begin{lema}\label{SumFormLeq}
Let $R$ be a finite ring and $ \alpha, x  \in R$.  If $ (x)_{\ell} = (x)_{\ell} \cap \langle \alpha \rangle^{\perp} $ then $$C_{\alpha}(x) = |[x]_{\ell}|.$$
\end{lema}
\begin{proof}
Assume that $ (x)_{\ell} = (x)_{\ell} \cap \langle \alpha \rangle^{\perp} $.  We have
\begin{align}
C_{\alpha}(x) &= \sum_{ \substack{ I \in \mathcal{I}(R) \\ I \subseteq (x)_{\ell} } }  |I|~  \mu_R( I , (x)_{\ell} ) \nonumber\\ 
&= |[x]_{\ell}|.\nonumber
\end{align}
Here the first equality follows from  Equation~(\ref{SubFormulaEq1}) and  the second equality follows from Lemma~\ref{GroupSizeEqEquiLema}.
\end{proof}

%%%%%%%%%%%%%%%%%%%%%%%%%%%%%%%%%%%%%%%%%%%%%%%%%%%%%%%%%%%%
\noindent Combining Lemma \ref{NewLemmaqqq} and Lemma \ref{SumFormLeq}, we obtain the following Theorem.
\begin{theorem}\label{MainResultProof} Let $R$ be a finite ring and $\alpha, x \in R$. Then 
\begin{align}
C_{\alpha}(x) = \left\{ \begin{array}{cl}
		  |[x]_{\ell}|  &  \mbox{if } K=(x)_{\ell}\\ 
	    \sum\limits_{ \substack{ E\subseteq \mathcal{M}_R((x)_{\ell}) \\ \bigcap\limits_{M \in E}M \subseteq  K } }  |\bigcap\limits_{M \in E}M|~  (-1)^{|E|}   & \mbox{if } K =\bigcap\limits_{M \in E } M  \mbox{ for some } E \subseteq \mathcal{M}_R((x)_{\ell}) \\
        0 &  \mbox{otherwise}  \\ 
	\end{array}\right. \nonumber
\end{align} 
where $K$ is the largest left ideal of $R$ contained in $(x)_{\ell} \cap \langle \alpha \rangle^{\perp}$. 
\end{theorem}

\begin{proof}
Assume that $K$ is the largest left ideal of $R$ contained in $(x)_{\ell} \cap \langle \alpha \rangle^{\perp}$. We have the following three cases:

  \noindent \textbf{Case 1:} Assume that $K=(x)_{\ell}$. The proof follows from Lemma~\ref{SumFormLeq}.

  \noindent\textbf{Case 2:} Assume that $K$ can be express as an intersection of some maximal $R$-left ideals of $ (x)_{\ell}$. The proof follows from Lemma~\ref{NewLemmaqqq}.

 \noindent\textbf{Case 3:} Assume that $K \neq (x)_{\ell}$ and $K$ can not be expressed as an intersection of some maximal $R$-left ideals of $ (x)_{\ell}$. From Equation~(\ref{NewEqUniqeKaddeq}), we have $$C_{\alpha}( x ) = \sum_{ \substack{ I \in \mathcal{I}(R) \\ I \subseteq  K } }  |I|~  \mu_R( I , (x)_{\ell} ).$$ By Lemma~\ref{MuZeroCond}, we have $\mu_R( I , (x)_{\ell} ) = 0$ for all $I \in \mathcal{I}(R)$ and $I \subseteq K$. Hence $C_{\alpha}(x)=0$.

\end{proof}

%%%%%%%%%%%%%%%%%%%%%%%%%%%%%%%%%%%%%%%%%%%%%%%%%%%%%%%%%%%%%%%%%%%%%%%%%%%%%%%%%%%%%%%%%%%%%%%%%%%%%%%%%%%%%%%%%%%%%%%%%%%%%%%%%%%%%%%%%%%%%%%%%%%%%%%%%%%%%%%%
%%%%%%%%%%%%%%%%%%%%%%%%%%%%%%%%%%%%%%%%%%%%%%%%%%%%%%%%%%%%%%%%%%%%%%%%%%%%%%%%
%%%%%%%%%%%%%%%%%%%%%%%%%%%%%%%%%%%%%%%%%%%%%%%%%%%%%%%%%%%%%%%%%%%%%%%%%%%%%%%%

The main result of this paper is the following.

\begin{theorem}\label{MainResultProof00} Let $R$ be a finite ring and $\alpha, x \in R$. If $\mathcal{M}_{R}((x)_{\ell})$ is minimal subset of itself, then
\begin{align}
C_{\alpha}(x) =\mu_R( K , (x)_{\ell} ) ~~ \frac{|[x]_{\ell}|}{\varphi_R( K , (x)_{\ell} )}. \nonumber
\end{align} 
where $K$ is the largest left ideal of $R$ contained in $(x)_{\ell} \cap \langle \alpha \rangle^{\perp}$. 
\end{theorem}

\begin{proof}
Assume that $\mathcal{M}_{R}((x)_{\ell})$ is minimal subset of itself and $K$ is the largest left ideal of $R$ contained in $(x)_{\ell} \cap \langle \alpha \rangle^{\perp}$. Let $\mathcal{M}_1 
 = \mathcal{M}_R( (x)_{\ell} ,~ K )$, and $\mathcal{M}_2 = \mathcal{M}_R((x)_{\ell}) \setminus \mathcal{M}_R((x)_{\ell},~ K )$. We have the following three cases:
 
 \noindent \textbf{Case 1:} Assume that $K=(x)_{\ell}$. Then $\mu_R( K , (x)_{\ell} ) = 1$ and $\varphi_R( K , (x)_{\ell} ) = 1$, and so the proof follows from Theorem~\ref{MainResultProof}.

\noindent\textbf{Case 2:} Assume that $K$ can be expressed as the intersection of some maximal $R$-left ideals of $(x)_{\ell}$. Since $\mathcal{M}_R((x)_{\ell})$ is a minimal subset of itself, it follows that $\mathcal{M}_1$ is the unique subset of $\mathcal{M}_R((x)_{\ell})$ satisfying  
\[
K = \bigcap_{M \in \mathcal{M}_1} M .
\]

Let $I \in \mathcal{I}(R)$ with $I \subseteq K$ and $\mu_R(I,(x)_{\ell}) \neq 0$. Then $I$ can also be expressed as the intersection of some maximal $R$-left ideals of $(x)_{\ell}$. Thus, there exists a unique subset $E_1 \subseteq \mathcal{M}_R((x)_{\ell})$ such that  
\[
I = \bigcap_{M \in E_1} M .
\]  
This implies that  
\[
I = I \cap K = \bigcap_{M \in E_1 \cup \mathcal{M}_1} M .
\]  
Note that $E_1 \cap \mathcal{M}_1$ may or may not be empty. Since $\mathcal{M}_R((x)_{\ell})$ is minimal in itself, it follows that there exists a unique subset $E \subseteq \mathcal{M}_R((x)_{\ell})$ such that  
\[
E \cap \mathcal{M}_1 = \emptyset 
\quad \text{and} \quad
E \cup \mathcal{M}_1 = E_1 \cup \mathcal{M}_1 .
\]  
Hence,  
\[
I = \bigcap_{M \in \mathcal{M}_1} M \ \cap \ \bigcap_{M \in E} M .
\]

\noindent Moreover,  
\[
\mu_R(I,(x)_{\ell}) = (-1)^{|\mathcal{M}_1|} \, (-1)^{|E|} 
= \mu_R(K,(x)_{\ell}) \, (-1)^{|E|}.
\]  
Therefore, for each $I$ with $\mu_R(I,(x)_{\ell}) \neq 0$ occurring in the summation of Equation~\eqref{NewEqUniqeKaddeq}, we have
\begin{align}
I = \bigcap_{M\in \mathcal{M}_1}M \bigcap\limits_{ M \in E } M \mbox{ and } \mu_R(I,(x)_{\ell}) =  \mu_G( K , (x)_{\ell} ) (-1)^{|E|} \label{NewEqFinResult11}
\end{align}
for some $E \subseteq \mathcal{M}_2$. Equation~(\ref{NewEqUniqeKaddeq}) implies that
\begin{align}
C_{\alpha}( x ) &= \sum_{ \substack{ I \in \mathcal{I}(R) \\ I \subseteq K  } }  |I|~  \mu_R( I , (x)_{\ell} ) \nonumber\\
&= \mu_R( K , (x)_{\ell} ) \sum_{\substack{ E \subseteq \mathcal{M}_2 }} | \bigcap_{M\in \mathcal{M}_1}M  \bigcap\limits_{M \in E}M|  ~ (-1)^{|E|} \nonumber\\
&= \mu_R( K , (x)_{\ell} ) ~~ |(x)_{\ell}| ~~ \sum_{\substack{ E \subseteq \mathcal{M}_2 }} \frac{| \bigcap\limits_{M \in \mathcal{M}_1}M  \bigcap\limits_{M \in E}M|}{|(x)_{\ell}|}  ~ (-1)^{|E|} \nonumber\\
&=\mu_R( K , (x)_{\ell} ) ~ |(x)_{\ell}| ~ \prod_{M \in \mathcal{M}_1} \frac{|M|}{|(x)_{\ell}|} \left[ 1+  \sum_{\substack{ E \subseteq \mathcal{M}_2 \\ E \neq \emptyset } } \frac{| \bigcap\limits_{M \in E} M|}{|(x)_{\ell}|}   ~  (-1)^{|E|} \right].
\label{SubFormulaEq1135}
\end{align}
Here the second equality hold from Equation (\ref{NewEqFinResult11}) and the fourth equality follows from Part $(ii)$ of Lemma~\ref{LemmaMaximEquality00}. By Lemma~\ref{NewLemaForMainThm}, we have 
\begin{align}
    \varphi_R( K , (x)_{\ell} ) = \sum_{E \subseteq \mathcal{M}_1 } (-1)^{|\mathcal{M}_1|-|E|} \prod_{M \in E} \frac{|(x)_{\ell}|}{|M|}.
    \label{SubFormulaEq1136}
\end{align}
\noindent Equation~(\ref{SubFormulaEq1135}) and Equation~(\ref{SubFormulaEq1136}) implies that
\begin{align}
 \varphi_R( K , (x)_{\ell} ) ~C_{\alpha}( x ) &= \mu_R( K , (x)_{\ell} ) ~ |(x)_{\ell}| ~ \left[ \sum_{E \subseteq \mathcal{M}_1 }  (-1)^{|\mathcal{M}_1|-|E|} \prod_{M \in E} \frac{|(x)_{\ell}|}{|M|} \prod_{M \in \mathcal{M}_1} \frac{|M|}{|(x)_{\ell}|} \right] \nonumber\\
&~~~~ \left[ 1 +   \sum_{\substack{ E \subseteq \mathcal{M}_2 \\ E \neq \emptyset } }  \frac{| \bigcap\limits_{M \in E} M|}{|(x)_{\ell}|}   ~  (-1)^{|E|} \right] \nonumber\\
&=\mu_R( K , (x)_{\ell} ) ~|(x)_{\ell}| ~  \left[ 1 +  \sum_{ \substack{E \subseteq \mathcal{M}_1  \\ E \neq \emptyset }}  (-1)^{|E|} \prod_{M \in E} \frac{|M|}{|(x)_{\ell}|}  \right] ~\nonumber\\
&~~~~ \left[ 1 +  \sum_{\substack{ E \subseteq \mathcal{M}_2 \\ E \neq \emptyset } } \frac{| \bigcap\limits_{M \in E} M|}{|(x)_{\ell}|}    ~  (-1)^{|E|} \right] \nonumber\\
&=\mu_R( K , (x)_{\ell} )  ~ |(x)_{\ell}| ~\left[ 1 +  \sum_{\substack{ E \subseteq \mathcal{M}_1 \\ E \neq \emptyset }} \frac{| \bigcap\limits_{M \in E} M|}{|(x)_{\ell}|} (-1)^{|E|}  \right] ~\nonumber\\
&~~~~ \left[ 1+  \sum_{\substack{ E \subseteq \mathcal{M}_2 \\ E \neq \emptyset } } \frac{| \bigcap\limits_{M \in E} M|}{|(x)_{\ell}|}    ~  (-1)^{|E|} \right] \nonumber\\
&=\mu_R( K , (x)_{\ell} ) ~ |(x)_{\ell}| ~ \left[  1 + \sum_{\substack{ E \subseteq \mathcal{M}_1  \cup  \mathcal{M}_2 \\ E \neq \emptyset }}  \frac{| \bigcap\limits_{M \in E} M|}{|(x)_{\ell}|}  ~  (-1)^{|E|} \right] \nonumber\\
&=\mu_R( K , (x)_{\ell} ) ~ \left[  |(x)_{\ell}| + \sum_{\substack{ E \subseteq \mathcal{M}_1  \cup  \mathcal{M}_2 \\ E \neq \emptyset }}  | \bigcap\limits_{M \in E} M| ~  (-1)^{|E|} \right]\nonumber\\
&=\mu_R( K , (x)_{\ell} ) ~ \sum_{ \substack{ I \in \mathcal{I}(R) \\ I \subseteq (x)_{\ell}  } }  |I|~  \mu_R( I , (x)_{\ell} ) \nonumber\\
&= \mu_R( K , (x)_{\ell} ) ~|[x]_{\ell}|.
\label{SubFormulaEq1137}
\end{align}
Here the fourth equality holds from Part $(ii)$ of Lemma~\ref{LemmaMaximEquality00} and the last equality follows from Lemma~\ref{GroupSizeEqEquiLema}. Now, the proof follows from Equation~(\ref{SubFormulaEq1137}). 

\noindent\textbf{Case 3:} Assume that $K \neq (x)_{\ell}$ and $K$ can not be expressed as an intersection of some maximal $R$-left ideals of $ (x)_{\ell}$. Then $\mu_R( K , (x)_{\ell} ) = 0$, and so the proof follows from Theorem~\ref{MainResultProof}. 

\end{proof}

%%%%%%%%%%%%%%%%%%%%%%%%%%%%%%%%%%%%%%%%%%%%%%%%%%%%%%%%%%%%%%%%%%%%%%%%%%%%%%%%%%%%%%%%%%%%%%%%%%%%%%%%%%%%%%%%%%%%%%%%%%%%%%%%%%%%%%%%%%%%%%%%%%%%%%%%%%%%%%%%
%%%%%%%%%%%%%%%%%%%%%%%%%%%%%%%%%%%%%%%%%%%%%%%%%%%%%%%%%%%%%%%%%%%%%%%%%%%%%%%%
%%%%%%%%%%%%%%%%%%%%%%%%%%%%%%%%%%%%%%%%%%%%%%%%%%%%%%%%%%%%%%%%%%%%%%%%%%%%%%%%

If $R$ is a commutative ring with unity, then left ideal and two-sided ideal are the same. Therefore, $\mathcal{M}_{R}((x)_{\ell})$ will always be a minimal subset of itself, owing to the Chinese Remainder Theorem and Lemma~\ref{LemmaMaximEquality00}. We have the following result for commutative rings.

\begin{corollary}\label{MainResultProofCorollary00} Let $R$ be a commutative finite ring and $\alpha, x \in R$. Then
\begin{align}
C_{\alpha}(x) =\mu_R( K , (x)_{\ell} ) ~~ \frac{|[x]_{\ell}|}{\varphi_R( K , (x)_{\ell} )}. \nonumber
\end{align} 
where $K$ is the largest left ideal of $R$ contained in $(x)_{\ell} \cap \langle \alpha \rangle^{\perp}$. 
\end{corollary}

%%%%%%%%%%%%%%%%%%%%%%%%%%%%%%%%%%%%%%%%%%%%%%%%%%%%%%%%%%%%%%%%%%%%%%%%%%%%%%%%%%%%%%%%%%%%%%%%%%%%%%%%%%%%%%%%%%%%%%%%%%%%%%%%%%%%%%%%%%%%%%%%%%%%%%%%%%%%%%%%

%%%%%%%%%%%%%%%%%%%%%%%%%%%%%%%%%%%%%%%%%%%%%%%%%%%%%%%%%%%%%%%%%%%%%%%%%%%%%%%%
%%%%%%%%%%%%%%%%%%%%%%%%%%%%%%%%%%%%%%%%%%%%%%%%%%%%%%%%%%%%%%%%%%%%%%%%%%%%%%%%

\section{Conclusion}
In this section, we give a idea to present the main results of this paper in terms of right ideals and two-sided ideals of $R$. For any $x \in R$, we denote by $(x)_{r}$ and $(x)_{t}$ the right ideal and two-sided ideal of $R$ generated by $x$, respectively. Define 
\begin{align*} 
[x]_{r} &:= \{ y \in R: (y)_{r} = (x)_{r} \},\\ 
[x]_{t} &:= \{ y \in R: (y)_{t} = (x)_{t} \}.
\end{align*} It is important fact that if $x\in R^{\times}$ then $[x]_{r}=R^{\times}$, and $R^{\times} \subseteq [x]_{t}$. However, the set $[x]_{t}$ may or may not be equal to $R^{\times}$.

Define the character sums $C_{\alpha}^{r}(x)$ and $C_{\alpha}^{t}(x)$ corresponding to the right ideal and the two-sided ideal, respectively, by replacing the set  $[x]_{\ell}$ with $[x]_{r}$ and $[x]_{t}$ in the definition of $C_{\alpha}(x)$. That is, for any $\alpha , x \in R$, we have
\begin{align}
C_{\alpha}^r(x)= \sum_{s \in [x]_{r} } \psi_{\alpha}(s) \mbox{ and } C_{\alpha}^t(x)= \sum_{s \in [x]_{t} } \psi_{\alpha}(s).\nonumber
\end{align} 

A similar formula can also be derived for $C_{\alpha}^{r}(x)$ and $C_{\alpha}^{t}(x)$ by defining the M$\ddot{\text{o}}$bius function for right and two-sided ideals, along with the corresponding notations for left and two-sided ideals. The formula for $C_{\alpha}^r(x)$ will resemble that of $C_{\alpha}(x)$, with analogous notations for right ideals. However, the case of the two-sided ideal  is simpler, as the Chinese Remainder Theorem holds for two-sided ideals. Define the following notations: $\mathcal{I}^t(R)$, $R$-two-sided ideal, maximal $R$-two-sided ideal, $\mathcal{M}^t_R(J)$, $\mathcal{M}_R^t(J,I)$, $\mu_R^t( I,J)$, $\varphi_R^t( I , J)$ by replacing the word ``left'' in $\mathcal{I}(R)$, $R$-left ideal, maximal $R$-left ideal, $\mathcal{M}_R(J)$, $\mathcal{M}_R(J,I)$, $\mu_R( I,J)$, $\varphi_R( I , J)$ with the word ``two-sided'', respectively. The formula for the character sum $C_{\alpha}^{t}(x)$ is given in the following result.

\begin{theorem}\label{MainResultProof02} Let $R$ be a finite ring and $\alpha , x \in R$. Then
\begin{align}
C_{\alpha}^t(x) =  \mu_R^t( K , (x)_{t} ) ~~ \frac{|[x]_{t}|}{\varphi_R^t( K , (x)_{t} )} \nonumber
\end{align} 
where $K$ is the largest two-sided ideal of $R$ contained in $(x)_{t} \cap \langle \alpha \rangle^{\perp}$. 
\end{theorem}

If we compare Theorem~\ref{MainResultProof00} with Theorem~\ref{MainResultProof02}, we observe that in Theorem~\ref{MainResultProof00}, we assume that $\mathcal{M}_{R}((x)_{\ell})$ is a minimal subset of itself. However, such a condition is not required in Theorem~\ref{MainResultProof02} because it is already satisfied due to the Chinese Remainder Theorem.

\section*{Acknowledgements} The first author is supported by Senior Research Fellowship from CSIR, Government of India (File No. 09/1020(15619)/2022-EMR-I).

\noindent\textbf{Data Availability} No data was gathered or used in this paper, so a ``data availability statement'' is not applicable.

\noindent\textbf{Conflict of interest} The author states that there is no conflict of interest.

%%%%%%%%%%%%%%%%%%%%%%%%%%%%%%%%%%%%%%%%%%%%%%%%%%%%%%%%%%%%%%%%%%%%%%%%%%%%%%%%%%%%%%

%\bibliographystyle{plain}
%\bibliography{SeprateBibFile}

\end{document}